# A Small-Gain Theorem with Applications to Input/Output Systems, Incremental Stability, Detectability, and Interconnections

Brian Ingalls and Eduardo D. Sontag*
Dept. of Mathematics, Rutgers University, NJ
{ingalls,sontag}@math.rutgers.edu

November 12, 2018

**Abstract**

A general ISS-type small-gain result is presented. It specializes to a small-gain theorem for ISS operators, and it also recovers the classical statement for ISS systems in state-space form. In addition, we highlight applications to incrementally stable systems, detectable systems, and to interconnections of stable systems.

## 1 Introduction

A common feature of control analysis is the application of small-gain results. One is often faced with feedback laws or auxiliary systems which are connected to a plant. Small-gain theorems can be used to verify stability of the resulting "closed-loop" systems, under appropriate conditions. Small-gain theorems have a long history, starting with the seminal work of Zames, Sandberg, Safanov, and others (see e.g. [23, 24, 36, 37] as well as expositions in for example [4, 13, 30, 34, 35]).

Most of this classical work, though not all, applies to norm-based (linear) gains. Much recent work has focused on versions of small-gain theorems expressed in terms of "nonlinear gain functions," see [20], and in particular on results which incorporate explicit estimates on transient behavior, expressed in the ISS framework. The fundamental work along these lines was that of Jiang, Teel, and Praly in [5], which gave rise to an extended follow-up literature, see e.g. [6, 7, 8, 31, 33].

The present work is also concerned with ISS-type small-gain theorems. The purpose of this paper is to present a very general principle, expressed in abstract terms, which allows application to a wide variety of contexts, which we summarize next.

The input-to-state stability (ISS) property was introduced in [26]; it represents a natural notion of stability for nonlinear control systems. Systems which satisfy the ISS property exhibit trajectories which are asymptotically bounded by a nonlinear gain on the inputs. In addition, the ISS bound includes a transient term which allows for a bounded overshoot depending on the size of the initial condition. The theory of ISS systems now forms an integral part of several texts [9, 12, 13, 16, 17, 25], as well as expository and research articles, see e.g. [5, 14, 18, 21, 27].

Also introduced in [26] was an accompanying notion of input-to-output stability (IOS) which can be expressed either for systems with state space representation or for purely input/output systems. The result in this paper can be applied in either case – the statement is made in sufficient generality to allow applicability to several situations.

*Supported in part by US Air Force Grants F49620-98-1-0421 and F49620-01-1-0063

ISS notions can also be generalized to *incremental* stability properties, which characterize systems for which each trajectory converges asymptotically to every other trajectory. The incremental-ISS property was addressed in [2]. The small-gain result presented here can be applied immediately to incrementally stable systems.

Another useful generalization of the ISS and IOS properties is to notions of detectability. Introduced in [28, 29] the *input-output-to-state stability* (IOSS) property is a notion of zero-detectability for nonlinear systems (see also [15]). The recently introduced property of *input-measurement-to-error stability* (IMES) [11] is a further generalization to systems with two outputs. Applications of the small-gain theorem to such systems is outlined below.

Perhaps the most common use of small-gain results is in verifying the stability of an interconnection of stable systems, and this is often how the results are presented. Some small-gain theorems for interconnections satisfying ISS properties have appeared in the literature. Interconnections of IOS systems with finite-dimensional state space representations were addressed in [5]. ISS interconnections for time-delay systems were studied in [32]. A small-gain theorem for purely input/output systems satisfying the IOS property was announced in [10]. The result in this paper recovers each of these results when applied to the appropriate situation.

The contents of this paper are presented as follows. Notations and definitions are presented in Section 2. The small-gain Theorem is stated and proved in Section 3. Section 4 is devoted to applications of the main result to various notions of stability and to interconnections of stable systems. The proofs of some basic results are included in an appendix.

## 2 Notations and Definitions

To allow application to several situations, the main result will be stated in rather abstract terms.

Let $\mathcal{T}$ be a subgroup of $(\mathbb{R}, +)$ which will be referred to as a *time set*. In most applications, we would find $\mathcal{T} = \mathbb{R}$ (continuous time) or $\mathcal{T} = \mathbb{Z}$ (discrete time). Let $\mathcal{T}_+$ denote the subset of nonnegative elements of $\mathcal{T}$. By notational convention, in what follows all intervals are assumed to be restricted to $\mathcal{T}$. For instance

$$[a, b) = \{t \in \mathcal{T} : a \leq t < b\}$$

and similarly for open, closed or infinite intervals.

Fix a set $\mathcal{U}$ which we will call the set of inputs. We introduce the following notation.

**Definition 2.1** We say that a quadruple $(\tau, u(\cdot), x(\cdot), y(\cdot))$ is a *trajectory* if $0 < \tau \leq \infty$, $u(\cdot) \in \mathcal{U}$, and $x(\cdot)$, $y(\cdot)$ are functions from $[0, \tau)$ into $\mathbb{R}_{\geq 0}$. □

**Remark 2.2** In what follows, the familiar names of state, output, and input will be used for $x(\cdot)$, $y(\cdot)$, and $u(\cdot)$, in hopes of aiding intuition of application to IOS systems. (Moreover, the symbol $u(\cdot)$ will be used throughout for the "input", even though the elements of the set $\mathcal{U}$ need not be functions). The reader should keep in mind the more general definition, in which no underlying connection between $x(\cdot)$, $y(\cdot)$ and $u(\cdot)$ is presumed. Notably, there is no assumption of *causality* (i.e. non-causal inputs are allowed). Indeed, the "input" $u(\cdot)$ is allowed to be a function of the "state", a situation that will be considered in the applications to detectability notions described below. □

**Remark 2.3** The primary motivation for the definition of a trajectory is to generalize input-state-output triples for systems. For example, consider a system (with $\mathcal{T} = \mathbb{R}$) given by

$$\dot{x}(t) = f(x(t), u(t)), \qquad y(t) = h(x(t), u(t)), \tag{1}$$

where $x \in \mathbb{R}^n$, $f : \mathbb{R}^n \times \mathbb{R}^m \to \mathbb{R}^n$ is locally Lipschitz, and $h : \mathbb{R}^n \times \mathbb{R}^m \to \mathbb{R}^p$. Inputs are measurable locally essentially bounded functions from $\mathbb{R}_{\geq 0}$ into $\mathbb{R}^m$. In this case we will use $|\cdot|$ for Euclidean norm and $\|\cdot\|$ for (essential) supremum norm. Given any input $u$ and any initial condition $\xi$, we let $x(\cdot, \xi, u)$ denote the unique maximal solution of the initial value problem $\dot{x} = f(x, u), x(0) = \xi$, defined on a maximal subinterval $[0, T^{\max}_{\xi,u})$. We will denote the corresponding output as $y(\cdot, \xi, u)$, that is $y(\cdot, \xi, u) = h(x(\cdot, \xi, u), u(\cdot))$ on the domain of definition of the solution. In this case, one can identify the quadruples $(\tau, u(\cdot), |x(\cdot, \xi, u)|, |y(\cdot, \xi, u)|)$ as trajectories for any $\tau \in (0, T^{\max}_{\xi,u})$. Note that in this type of setting, the set $\mathcal{U}$ characterizes both the input-value set and the form of the admissible input functions. □

In what follows we will need expressions for the magnitude of the "input" and "output" signals over time intervals. Since we will need to compare the value of $y(t)$ to a measure of the signal $y(\cdot)$ over an interval containing $t$, we will use the supremum norm to measure outputs. For each $0 \leq a \leq b$ we will denote

$$\|y\|_{[a,b]} := \sup\{|y(t)| : t \in [a, b]\},$$

whenever $y(\cdot)$ is a real-valued function defined on the interval $[a, b]$.

When measuring the input signal, we can allow more freedom. Let $\mathcal{N}$ denote the set of functions $\mu : \mathcal{U} \times \mathbb{R}_{\geq 0} \times \mathbb{R}_{\geq 0} \to \mathbb{R}_{\geq 0}$ which satisfy the monotonicity condition

$$\mu(u(\cdot), b, c) \leq \mu(u(\cdot), a, d) \qquad \forall\, 0 \leq a \leq b \leq c \leq d < \infty$$

for any $u(\cdot) \in \mathcal{U}$. That is, for each fixed $u(\cdot)$, $\mu(u(\cdot), \cdot, \cdot)$ defines a monotone set function on the finite subintervals of $[0, \infty)$. We call such functions *input measures*. Given any input measure $\mu \in \mathcal{N}$ we will use the notation

$$\|u\|_{[a,b]} := \mu(u(\cdot), a, b).$$

The use of these general "measures" allows applications to different norms on input signals, e.g. supremum norms or various integral norms.

A function $\gamma : \mathbb{R}_{\geq 0} \to \mathbb{R}_{\geq 0}$ is *of class* $\mathcal{K}$ (or a "$\mathcal{K}$-function") if it is continuous, positive definite, and strictly increasing; and is *of class* $\mathcal{K}_\infty$ if in addition it is unbounded. A function $\rho : \mathbb{R}_{\geq 0} \to \mathbb{R}_{\geq 0}$ is of class $\mathcal{L}$ if it is continuous, decreasing, and tends to zero as its argument tends to $+\infty$. A function $\beta : \mathbb{R}_{\geq 0} \times \mathbb{R}_{\geq 0} \to \mathbb{R}_{\geq 0}$ is *of class* $\mathcal{KL}$ if for each fixed $t \geq 0$, $\beta(\cdot, t)$ is of class $\mathcal{K}$ and for each fixed $s \geq 0$, $\beta(s, \cdot)$ is of class $\mathcal{L}$.

**Definition 2.4** Let $\mu \in \mathcal{N}$ and $C \geq 0$. A set $S$ of trajectories is called $(\mu, C)$-$\mathcal{KL}$-*practical-IOS* if there exists $\beta \in \mathcal{KL}$ so that for each $(\tau, u(\cdot), x(\cdot), y(\cdot)) \in S$ and each $t_0 \in [0, \tau)$,

$$y(t) \leq \max\{\beta(x(t_0), t - t_0), \|u\|_{[t_0, t]}, C\} \qquad \forall t \in [t_0, \tau).$$

□

**Definition 2.5** Let $\mu \in \mathcal{N}$ and $C \geq 0$. A set $S$ of trajectories is called $(\mu, C)$-*practical-IOS* if the following two properties hold.

1. $(\mu, C)$-uniform practical stability: There exists a $\mathcal{K}_\infty$-function $\delta(\cdot)$ such that for each $\varepsilon > 0$ and each $(\tau, u(\cdot), x(\cdot), y(\cdot)) \in S$, if $t_0 \in [0, \tau)$ is such that $x(t_0) \leq \delta(\varepsilon)$, then

$$y(t) \leq \max\{\varepsilon, \|u\|_{[t_0, t]}, C\} \qquad \forall t \in [t_0, \tau).$$

2. $(\mu, C)$-uniform practical attractivity: For any $r > 0$, $\varepsilon > C$, there is a $T = T_{r,\varepsilon} > 0$ so that for each $(\tau, u(\cdot), x(\cdot), y(\cdot)) \in S$, if $t_0 \in [0, \tau)$ is such that $x(t_0) \leq r$, then

$$y(t) \leq \max\{\varepsilon, \|u\|_{[t_0, t]}\} \qquad \forall t \in [t_0 + T, \tau).$$

□

The following minor generalization of Proposition 2.5 in [19] will be used in the proof of the small-gain Theorem. The proof of this result is included in the appendix.

**Proposition 2.6** Let $\mu \in \mathcal{N}$ and $C \geq 0$.
(i): If a set $S$ of trajectories is $(\mu, C)$-$\mathcal{KL}$-practical-IOS, then it is $(\mu, C)$-practical-IOS.
(ii): If a set $S$ of trajectories is $(\mu, C)$-practical-IOS, then it is $(\mu, 3C)$-$\mathcal{KL}$-practical-IOS. □

# 3 Small-Gain Theorem

This small-gain result complements Lemma A.1 in [5] and Lemma 3.4 in [15].

**Theorem 1** *Let $\mu^a$ and $\mu^b$ be input measures. Suppose given a $\mathcal{KL}$-function $\beta$, a number $r_0 \geq 0$, a $\mathcal{K}$-function $\gamma$ for which $\gamma(r) < r$ if $r > r_0$, $\mathcal{K}$-functions $\sigma_1$, $\sigma_2$, $\sigma_3$, and a constant $d \geq 0$. If there is a constant $C \geq 0$ and a set $S$ of trajectories each of which satisfy the following:*

*(H1) for each $t_0 \in [0, \tau)$,*

$$y(t) \leq \max\{\beta(x(t_0), t - t_0), \gamma(\|y\|_{[t_0,t]}), \|u\|_{[t_0,t]}^a, C\} \qquad \forall t \in [t_0, \tau),$$

*(H2) for each $t_0 \in [0, \tau)$,*

$$x(t) \leq \max\{\sigma_1(x(t_0)), \sigma_2(t - t_0), \sigma_3(\|y\|_{[t_0,t]}), \|u\|_{[t_0,t]}^b, d\} \qquad \forall t \in [t_0, \tau);$$

*then there exists a $\mathcal{KL}$-function $\widetilde{\beta}$ so that each trajectory $(\tau, u(\cdot), x(\cdot), y(\cdot))$ in $S$ also satisfies, for each $t_0 \in [0, \tau)$,*

$$y(t) \leq \max\{\widetilde{\beta}(x(t_0), t - t_0), \|u\|_{[t_0,t]}^m, 3C, 3r_0\} \qquad \forall t \in [t_0, \tau),$$

*where $\|\cdot\|_{\mathcal{I}}^m = \max\{\beta(\sigma_3(\|\cdot\|_{\mathcal{I}}^a), 0), \beta(\|\cdot\|_{\mathcal{I}}^b, 0), \|\cdot\|_{\mathcal{I}}^a\}$ for each interval $\mathcal{I} \subseteq \mathcal{T}_+$.* ■

**Remark 3.1** It is clear from the statement that the Theorem can be applied to a number of situations. One case which may not be transparent is the application to systems where the outputs are only measurable functions (rather than "true" functions) on their domain (as is the case when $y(\cdot)$ is a function of $u(\cdot)$, and $u(\cdot)$ is only measurable). The Theorem can be applied in this case by replacing $y(\cdot)$ with a "true" function $y^*(\cdot)$ from the equivalence class represented by $y(\cdot)$ such that

$$\sup\{|y^*(t)| : a \leq t \leq b\} = \text{ess sup}\{|y(t)| : a \leq t \leq b\}$$

for all intervals $[a, b]$. This can always be done, e.g. by setting

$$y^*(t) := \begin{cases} \lim_{h \to 0} \frac{\int_t^{t+h} y(s)\, ds}{h} & \text{if the limit exists} \\ 0 & \text{otherwise.} \end{cases}$$

In this case $y^*(t) = y(t)$ almost everywhere, by Lebesgue's Differentiation Theorem. □

*Proof.* (Theorem 1)

Let a $\mathcal{KL}$-function $\beta$, a number $r_0 \geq 0$, a $\mathcal{K}$-function $\gamma$ so that $\gamma(r) < r$ for each $r > r_0$, $\mathcal{K}$-functions $\sigma_1$, $\sigma_2$, $\sigma_3$, and a number $d \geq 0$ be given. Let $C \geq 0$ be fixed, and let $S$ be a set of trajectories each of which satisfy hypotheses (H1) and (H2). We are required to show that there exists $\widetilde{\beta} \in \mathcal{KL}$ so that every trajectory in $S$ is $(\mu^m, 3\max\{C, r_0\})$-$\mathcal{KL}$-practical-IOS. We will use Proposition 2.6, by which we need only show that each trajectory in $S$ satisfies conditions (1) and (2) of the definition of $(\mu^m, \max\{C, r_0\})$-practical-IOS.

Let $\vartheta(s) := \beta(s, 0)$ for all $s \geq 0$. Let $\widehat{C} := \max\{C, r_0\}$. Denote the composition of $\gamma$ with itself $n$ times by $\gamma^n$ (and likewise $\gamma^0(r) = r$).

*Claim:* The following holds for each trajectory in $S$. If $t_0 \in [0, \tau)$ is such that $x(t_0) \leq r$, then

$$y(t) \leq \max\{\vartheta(r), \|u\|_{[t_0,t]}^a, \widehat{C}\} \qquad \forall t \in [t_0, \tau). \tag{2}$$

*Proof of Claim:* Suppose $t_0$ and $(\tau, u(\cdot), x(\cdot), y(\cdot)) \in S$ are such that $x(t_0) \leq r$. Fix any $t \in [t_0, \tau)$. From (H1), we find that for any $s \in [t_0, t]$,

$$y(s) \leq \max\{\beta(x(t_0), s - t_0), \gamma(\|y\|_{[t_0,s]}), \|u\|_{[t_0,s]}^a, C\}.$$

So
$$\|y\|_{[t_0,t]} \leq \max\{\vartheta(r), \gamma(\|y\|_{[t_0,t]}), \|u\|_{[t_0,t]}^a, C\}.$$

Since $\gamma(r) < r$ for all $r > r_0$, (2) follows.

We next define some sequences. Fix $r > 0$. Let $T_0^r = 0$. Set
$$M_0^r := r,$$
and let $T_1^r$ be any element of $\mathcal{T}_+$ so that $\beta(M_0^r, T_1^r) \leq \gamma(\vartheta(r))$. Define $\hat{\sigma}_1(s) := \max\{\sigma_1(s), \sigma_3(\vartheta(s))\}$ for all $s \geq 0$, and let $\hat{d} := \max\{\sigma_3(\widehat{C}), d\}$. We then make the recursive definitions for each integer $i \geq 1$:
$$M_i^r := \max\{\hat{\sigma}_1(r), \sigma_2(T_0^r + \cdots + T_i^r), \hat{d}\}$$
and $T_{i+1}^r$ is any element of $\mathcal{T}_+$ such that $\beta(M_i^r, T_{i+1}^r) \leq \gamma^{i+1}(\vartheta(r))$. Finally, for each $i \geq 0$, define $\widehat{T}_i^r := T_0^r + T_1^r + T_2^r + \cdots + T_i^r$.

*Claim:* The following holds for each trajectory in $S$. If $t_0 \in [0, \tau)$ is such that $x(t_0) \leq r$, then for each $i \geq 0$,
$$x(t_0 + \widehat{T}_i^r) \leq \max\{M_i^r, \sigma_3(\|u\|_{[t_0, t_0 + \widehat{T}_i^r]}^a), \|u\|_{[t_0, t_0 + \widehat{T}_i^r]}^b\}, \tag{3}$$
provided $t_0 + \widehat{T}_i^r < \tau$.

*Proof of Claim:* Suppose $t_0$ and $(\tau, u(\cdot), x(\cdot), y(\cdot)) \in S$ are such that $x(t_0) \leq r$. Hypothesis (H2) and (2) give
$$\begin{aligned}
x(t) &\leq \max\{\sigma_1(r), \sigma_2(t - t_0), \sigma_3(\vartheta(r)), \sigma_3(\|u\|_{[t_0,t]}^a), \sigma_3(\widehat{C}), \|u\|_{[t_0,t]}^b, d\} \\
&= \max\{\hat{\sigma}_1(r), \sigma_2(t - t_0), \sigma_3(\|u\|_{[t_0,t]}^a), \|u\|_{[t_0,t]}^b, \hat{d}\}
\end{aligned}$$
for all $t \in [t_0, \tau)$. The claim follows from the definitions of $M_i^r$ and $\widehat{T}_i^r$.

*Condition (2): $(\mu^m, \widehat{C})$-uniform practical attractivity:* We next show that each trajectory of $S$ is $(\mu^m, \widehat{C})$-uniform practically attractive. We will first verify a decrease statement over the times $\widehat{T}_i^r$.

*Claim:* The following holds for each trajectory in $S$. If a time $t_0 \in [0, \tau)$ is such that $x(t_0) \leq r$, then, for each $i \geq 0$,
$$y(t) \leq \max\{\gamma^i(\vartheta(r)), \|u\|_{[t_0,t]}^m, \widehat{C}\} \tag{4}$$
for all $t \in [t_0 + \widehat{T}_i^r, \tau)$.

We will prove this by induction on the index $i$.

*Proof:* Suppose $t_0 \in [0, \tau)$ and $x(\cdot) \in S$ are such that $x(t_0) \leq r$. We have already shown (2), which gives (4) for $i = 0$.

Fix any $i \geq 1$, and suppose
$$y(t) \leq \max\{\gamma^{i-1}(\vartheta(r)), \|u\|_{[t_0,t]}^m, \widehat{C}\} \tag{5}$$
for all $t \in [t_0 + \widehat{T}_{i-1}^r, \tau)$.

Since (3) tells us that $x(t_0 + \widehat{T}_{i-1}^r) \leq \max\{M_{i-1}^r, \sigma_3(\|u\|_{[t_0, t_0 + \widehat{T}_{i-1}^r]}^a), \|u\|_{[t_0, t_0 + \widehat{T}_{i-1}^r]}^b\}$, condition (H1) gives that for each $t \in [t_0 + \widehat{T}_i^r, \tau)$,
$$\begin{aligned}
y(t) &\leq \max\{\beta(x(t_0 + \widehat{T}_{i-1}^r), t - (t_0 + \widehat{T}_{i-1}^r)), \gamma(\|y\|_{[t_0 + \widehat{T}_{i-1}^r, t]}), \|u\|_{[t_0 + \widehat{T}_{i-1}^r, t]}^a, C\} \\
&\leq \max\{\beta(M_{i-1}^r, t - (t_0 + \widehat{T}_{i-1}^r)), \beta(\sigma_3(\|u\|_{[t_0, t_0 + \widehat{T}_{i-1}^r]}^a), 0), \beta(\|u\|_{[t_0, t_0 + \widehat{T}_{i-1}^r]}^b, 0), \\
&\quad \gamma(\|y\|_{[t_0 + \widehat{T}_{i-1}^r, t]}), \|u\|_{[t_0 + \widehat{T}_{i-1}^r, t]}^a, C\} \\
&\leq \max\{\gamma^i(\vartheta(r)), \gamma(\|y\|_{[t_0 + \widehat{T}_{i-1}^r, t]}), \|u\|_{[t_0, t]}^m, C\}. \tag{6}
\end{aligned}$$

Now, take any $t \in [t_0 + \widehat{T}_i^r, \tau)$. Consider two cases.

Case (i): If
$$\gamma(\|y\|_{[t_0+\widehat{T}_{i-1}^r, t]}) < \max\{\gamma^i(\vartheta(r)), \|u\|_{[t_0,t]}^m, C\},$$
then (6) immediately gives
$$\begin{aligned} y(t) &\leq \max\{\gamma^i(\vartheta(r)), \|u\|_{[t_0,t]}^m, C\} \\ &\leq \max\{\gamma^i(\vartheta(r)), \|u\|_{[t_0,t]}^m, \widehat{C}\}. \end{aligned}$$

Case (ii): If
$$\gamma(\|y\|_{[t_0+\widehat{T}_{i-1}^r, t]}) \geq \max\{\gamma^i(\vartheta(r)), \|u\|_{[t_0,t]}^m, C\},$$
then (6) and the induction hypothesis (5) give
$$\begin{aligned} y(t) &\leq \gamma(\|y\|_{[t_0+\widehat{T}_{i-1}^r, t]}) \\ &\leq \gamma(\max\{\gamma^{i-1}(\vartheta(r)), \|u\|_{[t_0,t]}^m, \widehat{C}\}) \\ &\leq \max\{\gamma^i(\vartheta(r)), \gamma(\|u\|_{[t_0,t]}^m), \gamma(\widehat{C})\} \\ &\leq \max\{\gamma^i(\vartheta(r)), \|u\|_{[t_0,t]}^m, \widehat{C}, r_0\} \\ &= \max\{\gamma^i(\vartheta(r)), \|u\|_{[t_0,t]}^m, \widehat{C}\}. \end{aligned}$$

The claim follows by induction.

Now, for each fixed $r > 0$, we consider $\{T_r^i\}_{i=0}^\infty$ and $\{\widehat{T}_r^i\}_{i=0}^\infty$ as defined above. For each $\varepsilon > \widehat{C}$, we set $T_{r,\varepsilon} = \widehat{T}_i^r$ where $i$ is the smallest index such that $\gamma^i(\vartheta(r)) < \varepsilon$. Such an index always exists since for each $r > r_0$, the sequence $\{\gamma^n(r)\}_{n=0}^\infty$ is decreasing as long as its elements are greater than $r_0$, and $\limsup_{n \to \infty} \gamma^n(r) \leq r_0$ (since each $r \leq r_0$ must have $\gamma(r) \leq r_0$, which follows from the fact that $\gamma$ is increasing and that, by continuity, $\gamma(r_0) \leq r_0$).

To complete the proof, we next show the stability property.

*Condition (1): $(\mu^m, \widehat{C})$-uniform practical stability:* We will show that each trajectory in $S$ satisfies the $(\mu^m, \widehat{C})$-uniform practical stability property. Recall that (2) gives, in particular, that for any $r > 0$, any $t_0 \in [0, \tau)$, and any $(\tau, u(\cdot), x(\cdot), y(\cdot)) \in S$ so that $x(t_0) \leq r$,
$$y(t) \leq \max\{\vartheta(r), \|u\|_{[t_0,t]}^m, \widehat{C}\} \qquad \forall t \in [t_0, \tau).$$

Thus we can take $\delta(\cdot) \in \mathcal{K}_\infty$ so that
$$\delta(\varepsilon) \leq \vartheta^{-1}(\varepsilon)$$
for $\varepsilon \in [0, \sup_{s \geq 0} \vartheta(s))$, with which each trajectory $(\tau, u(\cdot), x(\cdot), y(\cdot)) \in S$ satisfies the $(\mu^m, \widehat{C})$-uniform practical stability property.

Finally, we invoke Proposition 2.6 to conclude that the set $S$ is $(\mu^m, 3\widehat{C})$-$\mathcal{KL}$-practical-IOS with some $\widetilde{\beta} \in \mathcal{KL}$. This is the desired result. ∎

**Remark 3.2** Typically, bounds such as (H1) and (H2) involve either a maximum or a sum on the right-hand-side. The statement of (H1) as a maximum allowed for an efficient presentation of the small-gain result. However, it is commonly the case that such a bound naturally presents itself as a sum of terms. Theorem 1 can be applied in such a situation by rewriting the bound as a maximum. In this case, however, the hypothesis that $\gamma$ is a contraction must be strengthened to allow for a margin, since the gain must increase in size when the bound is rewritten.

We use the following elementary observation: if $\rho$ is a function of class $\mathcal{K}_\infty$, then, for each $a, b \geq 0$,
$$a + b \leq \max\{a + \rho(a), \rho^{-1}(b) + b\}.$$

This follows immediately by considering the two cases $b \leq \rho(a)$ and $b > \rho(a)$. □

**Corollary 3.3** Let $\mu^a$ and $\mu^b$ be input measures. Suppose given a $\mathcal{KL}$-function $\beta$, a number $r_0 \geq 0$, a $\mathcal{K}$-function $\gamma$ and a $\mathcal{K}_\infty$-function $\rho$ for which

$$\gamma(r) + \rho(\gamma(r)) < r \qquad \forall r > r_0,$$

$\mathcal{K}$-functions $\sigma_1$, $\sigma_2$, $\sigma_3$, and a constant $d \geq 0$. If there is a constant $C \geq 0$ and a set $S$ of trajectories each of which satisfy the following:

(H1') for each $t_0 \in [0, \tau)$,

$$y(t) \leq \beta(x(t_0), t - t_0) + \gamma(\|y\|_{[t_0, t]}) + \|u\|^a_{[t_0, t]} + C \qquad \forall t \in [t_0, \tau) \tag{7}$$

(H2) for each $t_0 \in [0, \tau)$,

$$x(t) \leq \max\{\sigma_1(x(t_0)), \sigma_2(t - t_0), \sigma_3(\|y\|_{[t_0, t]}), \|u\|^b_{[t_0, t]}, d\} \qquad \forall t \in [t_0, \tau);$$

then there exists a $\mathcal{KL}$-function $\widetilde{\beta}$ so that each trajectory $(\tau, u(\cdot), x(\cdot), y(\cdot))$ in $S$ also satisfies, for each $t_0 \in [0, \tau)$,

$$y(t) \leq \max\{\widetilde{\beta}(x(t_0), t - t_0), \|u\|^m_{[t_0, t]}, 3\alpha(C), 3r_0\} \qquad \forall t \in [t_0, \tau)$$

where $\alpha(s) := \max\{4\rho^{-1}(3s), 4s\}$ and $\|\cdot\|^m_\mathcal{I} = \max\{\alpha(\beta(\sigma_3(\alpha(\|\cdot\|^a_\mathcal{I}), 0))), \alpha(\beta(\|\cdot\|^b_\mathcal{I}, 0)), \alpha(\|\cdot\|^a_\mathcal{I})\}$ for each interval $\mathcal{I} \subseteq \mathcal{T}_+$. □

*Proof.* The bound (7) gives

$$\begin{aligned} y(t) &\leq \max\{\gamma(\|y\|_{[t_0,t]}) + \rho(\gamma(\|y\|_{[t_0,t]})), \rho^{-1}(\beta(x(t_0), t - t_0) + \|u\|^a_{[t_0,t]} + C) \\ &\qquad + \beta(x(t_0), t - t_0) + \|u\|^a_{[t_0,t]} + C\} \\ &\leq \max\{\alpha(\beta(x(t_0), t - t_0)), \gamma(\|y\|_{[t_0,t]}) + \rho(\gamma(\|y\|_{[t_0,t]})), \alpha(\|u\|^a_{[t_0,t]}), \alpha(C)\}. \end{aligned}$$

Theorem 1 can now be applied, since $\gamma + \rho \circ \gamma$ is a contraction. ∎

## 4 Applications

We describe several applications of Theorem 1.

### 4.1 Unboundedness Observable Systems

Recall the definition of system (1). The system is called *forward complete* if for every initial condition $\xi$ and every input $u$ the solution $x(\cdot, \xi, u)$ is defined on all of $\mathbb{R}_{\geq 0}$. A weaker property, introduced in [22], is the following.

**Definition 4.1** The system (1) is said to be *unboundedness observable* (UO) if every trajectory $x(\cdot, \xi, u)$ which has a finite maximal domain of definition $[0, T^{\max}_{\xi, u})$ satisfies

$$\limsup_{t \to T^{\max}_{\xi,u}} |y(t, \xi, u)| = \infty.$$

□

**Remark 4.2** As shown in Lemma 2.2 of [3], this definition of unboundedness observability can be expressed equivalently as a statement of the form (H2), which bounds the state in terms of the initial condition, time, and the output and input signals. This definition of UO is slightly weaker than the definition of UO given [5], which bounds the state in terms of only the initial condition, input, and output. □

When applying Theorem 1 to systems which satisfy the unboundedness observability condition, hypothesis (H2) is immediate, by Lemma 2.2 of [3]. Of course the same is true of application to forward complete systems or to systems where the output is the state (i.e. with $h$ the identity map), since such systems are always unboundedness observable.

## 4.2 Incremental Stability

Given a system as in (1), one can define *incremental stability* of the system as stability of the trajectories *to each other*. The following definition has appeared in the ISS framework (see also [29] for an earlier dual version for detectability).

**Definition 4.3** (cf. [2]) A system as in (1) is said to be *incrementally-ISS* if there exist $\beta \in \mathcal{KL}$ and $\gamma \in \mathcal{K}$ so that for each pair of initial conditions $\xi_1$, $\xi_2$, and each pair of inputs $u_1$, $u_2$, the trajectories satisfy

$$|x(t, \xi_1, u_1) - x(t, \xi_2, u_2)| \leq \beta(|\xi_1 - \xi_2|, t) + \gamma(\|u_1 - u_2\|_{[0,t]})$$

for all $t$ in the interval $[0, T^{\max}) := [0, \min\{T^{\max}_{\xi_1, u_1}, T^{\max}_{\xi_2, u_2}\})$. □

Theorem 1 yields a small-gain result for incrementally stable systems as follows. By thinking of a super-system consisting of two copies of (1), one can interpret a pair

$$(x_1(\cdot), x_2(\cdot)) := (x(\cdot, \xi_1, u_1), x(\cdot, \xi_2, u_2))$$

as a single trajectory which corresponds to the initial condition $(\xi_1, \xi_2)$ and input $(u_1, u_2)$. Theorem 1 may be applied to this super-system by choosing trajectories of the form $(T^{\max}, u_1 - u_2, x_1(\cdot) - x_2(\cdot), x_1(\cdot) - x_2(\cdot))$. In this case, hypothesis (H2) is always satisfied, since the "state" and the "output" coincide.

## 4.3 Detectability Notions

Several notions of detectability (or more precisely zero-detectability) have been formulated within the ISS framework. A basic definition is the following, see [15].

**Definition 4.4** The system (1) is said to be *input-output-to-state stable* (IOSS) if there exist $\beta \in \mathcal{KL}$ and $\gamma_1, \gamma_2 \in \mathcal{K}$ so that for each initial condition $\xi$, and each input $u$ the trajectories satisfy

$$|x(t, \xi, u)| \leq \beta(|\xi|, t) + \gamma_1(\|y\|_{[0,t]}) + \gamma_2(\|u\|_{[0,t]})$$

for all $t$ in the interval $[0, T^{\max}_{\xi, u})$. □

Theorem 1 can be applied in this case if one chooses trajectories of the form $(T^{\max}_{\xi,u}, x(\cdot), x(\cdot), (y(\cdot), u(\cdot)))$, i.e. both the "state" and the "output" correspond to $x$, and the "input" corresponds to the pair $(y, u)$.

The IOSS property can be generalized further to provide a notion of *partial detectability* for systems with two outputs. Consider the augmentation of system (1) by the addition of a second output

$$w(t) = k(x(t), u(t)).$$

If the output $y$ corresponds to a measurement, while the output $w$ indicates an error which is to be regulated, one might be interested in characterizing the notion of partial detectability of $w$ through $y$. The recently introduced property of *input-measurement-to-error stability* (IMES) provides such a notion (see [11]). The definition is the same Definition 4.4 above, except the state $x(t, \xi, u)$ on the left-hand-side is replaced by the error $w(t, \xi, u)$. Again, Theorem 1 can be applied immediately, simply by choosing trajectories of the form $(T^{\max}_{\xi,u}, x(\cdot), w(\cdot), (y(\cdot), u(\cdot)))$

## 4.4 Input/output systems

When considering the analogy of an IOS statement such as (H1) to the case of purely input/output systems, it is natural to identify the "state" $x$ at time $t_0$ with the "input so far" at $t_0$ (e.g. $x(t_0) := \|u\|_{[0,t_0]}$). This was the procedure followed in the definition of IOS for input/output systems given in [26]. In such cases it is immediate that (H2) holds (with $\sigma_3$ the identity and $\sigma_1$, $\sigma_2$ and $d$ arbitrary).

A small gain result for such systems was presented in [10]. That result is a consequence of Theorem 1 as will be shown in Section 4.5.2.

## 4.5 Interconnections

Perhaps the most common application of small-gain results is to interconnections of systems. Indeed the small-gain results in [5] and [10] are stated in that form.

To apply Theorem 1 to such interconnections, one simply considers the interconnection as a single super-system, with "coupled" input and output pairs. The small-gain condition appears as the typical requirement that an appropriate composition of gains is a contraction.

A number of results on interconnections are direct corollaries of Theorem 1. For ease of reference we present the small-gain results in [5] and [10] and indicate how they follow from Theorem 1.

### 4.5.1 IOS – Finite Dimensional State Space Representation

We begin by showing how Theorem 1 can be applied to an interconnection of IOS systems which allow finite dimensional state space representations. This result first appeared as Theorem 2.1 of [5].

Expanding on system (1), suppose given an interconnected system of the form

$$\dot{x}_1 = f_1(x_1, y_2, u_1) \qquad y_1 = h_1(x_1, y_2, u_1), \tag{8}$$
$$\dot{x}_2 = f_2(x_2, y_1, u_2) \qquad y_2 = h_2(x_2, y_1, u_2), \tag{9}$$

whose outputs are described by

$$y_1 = h_1(x_1, h_2(x_2, y_1, u_1), u_1) \tag{10}$$
$$y_2 = h_2(x_2, h_1(x_1, y_2, u_2), u_2). \tag{11}$$

**Corollary 4.5** *Suppose that, when viewed as independent systems, (8) and (9) satisfy the hypotheses of Theorem 1 with $\sigma_2 = 0$. That is, there exist $\beta_1, \beta_2 \in \mathcal{KL}$, $\gamma_1^y, \gamma_2^y, \gamma_1^u, \gamma_2^u, \sigma_1^1, \sigma_2^1, \sigma_1^3, \sigma_2^3, \sigma_1^4, \sigma_2^4 \in \mathcal{K}$, and nonnegative constants $C_1$, $C_2$, $d_1$, $d_2$ so that for any initial conditions $(\xi_1, \xi_2)$ and inputs $(u_1, u_2)$,*

$$|y_1(t)| \leq \max\{\beta_1(|\xi_1|, t), \gamma_1^y(\|y_2\|_{[0,t]}), \gamma_1^u(\|u_1\|_{[0,t]}), C_1\} \tag{12}$$
$$|y_2(t)| \leq \max\{\beta_2(|\xi_2|, t), \gamma_2^y(\|y_1\|_{[0,t]}), \gamma_2^u(\|u_2\|_{[0,t]}), C_2\} \tag{13}$$

$$|x_1(t)| \leq \max\{\sigma_1^1(|\xi_1|), \sigma_1^3(\|y_1\|_{[0,t]}), \sigma_1^4(\|u_1\|_{[0,t]}), d_1\} \tag{14}$$
$$|x_2(t)| \leq \max\{\sigma_2^1(|\xi_2|), \sigma_2^3(\|y_2\|_{[0,t]}), \sigma_2^4(\|u_2\|_{[0,t]}), d_2\}, \tag{15}$$

*for all $t \in [0, T^{max}) := [0, T^{max}_{\xi_1, \xi_2 u_1, u_2})$. Suppose further that there exists $r_0 \geq 0$ so that*

$$\begin{array}{c} \gamma_1^y(\gamma_2^y(r)) < r \\ \gamma_2^y(\gamma_1^y(r)) < r \end{array} \qquad \forall r > r_0. \tag{16}$$

*Then the interconnection is practically IOS in the following sense. There exist $\beta \in \mathcal{KL}$, $\gamma \in \mathcal{K}$ and a nonnegative constant $C$ (where $C = 0$ when $C_1 = C_2 = d_1 = d_2 = r_0 = 0$) so that for any initial condition pair $(\xi_1, \xi_2)$ and any input pair $(u_1, u_2)$, if the pair $(y_1, y_2)$ satisfies (10) and (11) on the interval $[0, T^{max})$, then for any $t \in [0, T^{max})$,*

$$|(y_1(t), y_2(t))| \leq \max\{\beta(|(\xi_1, \xi_2)|, t), \gamma(\|(u_1, u_2)\|_{[0,t]}), C\},$$

*where $|(a, b)| := \max\{|a|, |b|\}$ and $\|(v, w)\| := \max\{\|v\|, \|w\|\}$.* □

**Remark 4.6** In the case where $r_0 = 0$, each of the inequalities in (16) implies the other; if $r_0 > 0$, the same is true with a possibly larger value of $r_0$.

The bounds (14) and (15), while stronger than the UO property defined in section 4.1, are equivalent to the definition of UO given in [5]. □

*Proof.* The result will follow by showing that the hypotheses of Theorem 1 are satisfied by the coupled system. We begin by noting that for any $t \in [0, T^{\max})$, (13) gives

$$\|y_2\|_{[0,t]} \leq \max\{\beta_2(|\xi_2|, 0), \gamma_2^y(\|y_1\|_{[0,t]}), \gamma_2^u(\|u_2\|_{[0,t]}), C_2\}, \tag{17}$$

and

$$\|y_2\|_{[t/2,t]} \leq \max\{\beta_2(|\xi_2|, t/2), \gamma_2^y(\|y_1\|_{[0,t]}), \gamma_2^u(\|u_2\|_{[0,t]}), C_2\}. \tag{18}$$

Together, (12) and (17) give, for any $t \in [0, T^{\max})$,

$$\begin{aligned}|y_1(t)| &\leq \max\{\beta_1(|\xi_1|, t), \gamma_1^y(\beta_2(|\xi_2|, 0)), \gamma_1^y(\gamma_2^y(\|y_1\|_{[0,t]})), \gamma_1^y(\gamma_2^u(\|u_2\|_{[0,t]})), \\ & \qquad \gamma_1^y(C_2), \gamma_1^u(\|u_1\|_{[0,t]}), C_1\},\end{aligned}$$

from which we conclude

$$\begin{aligned}\|y_1\|_{[0,t]} &\leq \max\{\beta_1(|\xi_1|, 0), \gamma_1^y(\beta_2(|\xi_2|, 0)), \gamma_1^y(\gamma_2^y(\|y_1\|_{[0,t]})), \\ & \qquad \gamma_1^y(\gamma_2^u(\|u_2\|_{[0,t]})), \gamma_1^u(\|u_1\|_{[0,t]}), \gamma_1^y(C_2), C_1\}.\end{aligned}$$

Then, from the fact that $\gamma_1^y(\gamma_2^y(r)) < r$ for all $r > r_0$, we have

$$\|y_1\|_{[0,t]} \leq \max\{\beta_1(|\xi_1|, 0), \gamma_1^y(\beta_2(|\xi_2|, 0)), \gamma_1^y(\gamma_2^u(\|u_2\|_{[0,t]})), \gamma_1^u(\|u_1\|_{[0,t]}), \gamma_1^y(C_2), C_1, r_0\} \tag{19}$$

for all $t \in [0, T^{\max})$. Combining (14) with (19), we have, for any $t \in [0, T^{\max})$,

$$\begin{aligned}|x_1(t)| &\leq \max\{\sigma_1^1(|\xi_1|), \sigma_1^3(\beta_1(|\xi_1|, 0)), \sigma_1^3(\gamma_1^y(\beta_2(|\xi_2|, 0))), \sigma_1^3(\gamma_1^y(\gamma_2^u(\|u_2\|_{[0,t]}))), \\ & \qquad \sigma_1^3(\gamma_1^u(\|u_1\|_{[0,t]})), \sigma_1^3(\gamma_1^y(C_2)), \sigma_1^3(C_1), \sigma_1^3(r_0), \sigma_1^4(\|u_1\|_{[0,t]}), d_1\}.\end{aligned} \tag{20}$$

For notational convenience, let $\widetilde{C}_1 = \max\{\sigma_1^3(\gamma_1^y(C_2)), \sigma_1^3(C_1), \sigma_1^3(r_0), d_1\}$, and define $\widetilde{C}_2$ analogously. Finally, we invoke (12) again. The bounds (18) and (20) allow us to conclude that for any $t \in [0, T^{\max})$,

$$\begin{aligned}|y_1(t)| &\leq \max\{\beta_1(|x_1(t/2)|, t/2), \gamma_1^y(\|y_2\|_{[t/2,t]}), \gamma_1^u(\|u_1\|_{[t/2,t]}), C_1\} \\ &\leq \max\{\beta_1(\sigma_1^1(|\xi_1|), t/2), \beta_1(\sigma_1^3(\beta_1(|\xi_1|, 0)), t/2), \beta_1(\sigma_1^3(\gamma_1^y(\beta_2(|\xi_2|, 0))), t/2), \\ & \qquad \beta_1(\sigma_1^3(\gamma_1^y(\gamma_2^u(\|u_2\|_{[0,t/2]}))), t/2), \beta_1(\sigma_1^3(\gamma_1^u(\|u_1\|_{[0,t/2]})), t/2), \\ & \qquad \beta_1(\sigma_1^4(\|u_1\|_{[0,t/2]}), t/2), \beta_1(\widetilde{C}, t/2), \gamma_1^y(\beta_2(|\xi_2|, t/2)), \gamma_1^y(\gamma_2^y(\|y_1\|_{[0,t]})), \\ & \qquad \gamma_1^y(\gamma_2^u(\|u_2\|_{[0,t]})), \gamma_1^y(C_2), \gamma_1^u(\|u_1\|_{[t/2,t]}), C_1\}.\end{aligned}$$

Since an analogous statement holds for $y_2$, we conclude that hypothesis (H1) of Theorem 1 is satisfied. That is

$$|(y_1(t), y_2(t))| \leq \max\{\beta(|(\xi_1, \xi_2)|, t), \gamma(\|(y_1, y_2)\|_{[0,t]}), \gamma^u(\|(u_1, u_2)\|_{[0,t]}), C\} \qquad \forall t \in [0, T^{\max}).$$

where

$$\begin{aligned}\beta(s,t) &:= \max\{\beta_1(\sigma_1^1(s), t/2), \beta_1(\sigma_1^3(\beta_1(s, 0)), t/2), \beta_1(\sigma_1^3(\gamma_1^y(\beta_2(s, 0))), t/2), \gamma_1^y(\beta_2(s, t/2)), \\ & \qquad \beta_2(\sigma_2^1(s), t/2), \beta_2(\sigma_2^3(\beta_2(s, 0)), t/2), \beta_2(\sigma_2^3(\gamma_2^y(\beta_1(s, 0))), t/2), \gamma_2^y(\beta_1(s, t/2))\} \\ \gamma(r) &:= \max\{\gamma_1^y(\gamma_2^y(r)), \gamma_2^y(\gamma_1^y(r))\} \\ \gamma^u(r) &:= \max\{\beta_1(\sigma_1^3(\gamma_1^y(\gamma_2^u(r))), 0), \beta_1(\sigma_1^3(\gamma_1^u(r)), 0), \beta_1(\sigma_1^4(r), 0), \gamma_1^y(\gamma_2^u(r)), \gamma_1^u(r), \\ & \qquad \beta_2(\sigma_2^3(\gamma_2^y(\gamma_1^u(r))), 0), \beta_2(\sigma_2^3(\gamma_2^u(r)), 0), \beta_2(\sigma_2^4(r), 0), \gamma_2^y(\gamma_1^u(r)), \gamma_2^u(r)\} \\ C &:= \max\{\beta_1(\widetilde{C}_1, 0), \gamma_1^y(C_2), C_1, \beta_2(\widetilde{C}_2, 0), \gamma_2^y(C_1), C_2\}.\end{aligned}$$

Hypothesis (H2) is immediate from the statements (14) and (15). The conclusion follows directly from an application of Theorem 1 to the interconnected system. ∎

### 4.5.2 IOS – Purely Input/Output Systems

As mentioned earlier, Theorem 1 can be applied directly to systems which are represented in purely input/output form. Although there are many meaningful ways in which to define i/o operators, we will focus on the definition used in [26]. The interconnection result presented in this section was first announced in [10].

Given any normed linear space $S$, the space of measurable locally essentially bounded maps $w : \mathbb{R} \to S$ for which there is some $t_0 \in \mathbb{R}$ such that

$$w(t) = 0 \quad \forall t < t_0.$$

will be denoted $L_0^\infty(S)$. The symbol $|\cdot|$ will be used for the norm on $S$ and $\|\cdot\|$ for the essential supremum norm of functions in $L_0^\infty(S)$. We interpret any product of normed linear spaces $S_1 \times S_2$ as a space with norm $|(s_1, s_2)| = \max\{|s_1|, |s_2|\}$.

**Definition 4.7** For any two normed linear spaces $W$ and $Y$, we will say that a map

$$F : L_0^\infty(W) \to L_0^\infty(Y)$$

is an *input/output (i/o) operator* if the following properties hold:

a. $F$ is causal (i.e. if $w(t) = v(t)$ for almost every $t \leq T$ then $F(w)(t) = F(v)(t)$ for almost every $t \leq T$); and

b. $F$ is shift invariant (i.e. $F(\sigma_t w) = \sigma_t F(w)$ for all $t \in \mathbb{R}$ and all $w$, where $\sigma_t$ is the shift operator defined by $\sigma_t w(s) = w(s - t)$).

We say that the i/o operator $F$ is *input/output stable* (IOS) if, in addition:

c. there exist $\beta \in \mathcal{KL}$ and $\gamma \in \mathcal{K}$ such that

$$|F(w)(t)| \leq \max\{\beta(\|w\|_{(-\infty,0]}, t), \gamma(\|w\|_{[0,t]})\} \qquad \text{a.e. } t \geq 0,$$

for all $w \in L_0^\infty(W)$.

$\square$

Suppose $U_1, U_2, Y_1, Y_2$ are normed linear spaces. Given a pair of i/o operators

$$F_1 : L_0^\infty(U_1 \times Y_2) \to L_0^\infty(Y_1)$$
$$F_2 : L_0^\infty(U_2 \times Y_1) \to L_0^\infty(Y_2)$$

we consider the interconnected system represented by the coupled set of equations

$$\begin{aligned} y_1(t) &= F_1(u_1, y_2)(t) \\ y_2(t) &= F_2(u_2, y_1)(t) \end{aligned} \qquad \text{a.e. } t \in \mathbb{R}. \tag{21}$$

**Corollary 4.8** *Suppose that the i/o operators $F_1$ and $F_2$ are IOS with bounds $\beta_1$, $\beta_2 \in \mathcal{KL}$ and $\gamma_1^u$, $\gamma_1^y$, $\gamma_2^u$, $\gamma_2^y \in \mathcal{K}$ so that*

$$|F_1(u_1, y_2)(t)| \leq \max\{\beta_1(\|(u_1, y_2)\|_{(-\infty,0]}, t), \gamma_1^u(\|u_1\|_{[0,t]}), \gamma_1^y(\|y_2\|_{[0,t]})\}$$

*and*

$$|F_2(u_2, y_1)(t)| \leq \max\{\beta_2(\|(u_2, y_1)\|_{(-\infty,0]}, t), \gamma_2^u(\|u_2\|_{[0,t]}), \gamma_2^y(\|y_1\|_{[0,t]})\}$$

are satisfied for all $u_i \in L_0^\infty(U_i)$, $y_i \in L_0^\infty(Y_i)$ ($i = 1, 2$), and almost every $t \geq 0$. Then, if the gains are such that

$$\begin{aligned} \text{either} \quad & \gamma_1^y(\gamma_2^y(s)) < s \quad \forall s > 0, \\ \text{or} \quad & \gamma_2^y(\gamma_1^y(s)) < s \quad \forall s > 0, \end{aligned}$$

then the interconnected system is IOS in the following sense. There exist $\beta \in \mathcal{KL}$ and $\gamma \in \mathcal{K}$ such that for any inputs $(u_1, u_2) \in L_0^\infty(U_1 \times U_2)$, we have, for any solution $(y_1, y_2)$ of (21) and for almost every $t \geq 0$,

$$|(y_1(t), y_2(t))| \leq \max\{\beta(\|(u_1, u_2)\|_{(-\infty, 0]}, t), \gamma(\|(u_1, u_2)\|_{[0, t]})\}.$$

$\square$

The proof is similar to that shown in the previous section. The trajectories considered are quadruples of the form

$$\begin{aligned} (\infty, \|(u_1, y_2)\|_{(-\infty, \cdot]}, y_1(\cdot), (u_1(\cdot), y_2(\cdot))) \\ (\infty, \|(u_2, y_1)\|_{(-\infty, \cdot]}, y_2(\cdot), (u_2(\cdot), y_1(\cdot))). \end{aligned}$$

As mentioned above, the role of the "state" at time $t$ is played by the "input so far" at time $t$, e.g. $x_1(t) = \|(u_1, y_2)\|_{(-\infty, t]}$.

# A  Appendix

Before giving the proof of Proposition 2.6, we state the following result, which is proved in [19], and is stated as Lemma 4.1 in [1] (that reference requires $\varphi$ to take nonnegative values, but this can always be assumed without loss of generality, simply replacing $\varphi$ by $\max\{\varphi, 0\}$).

**Proposition A.1** If a function $\varphi : \mathbb{R}_{\geq 0} \times \mathbb{R}_{\geq 0} \to \mathbb{R}$ satisfies

- for all $r > 0$, $\varepsilon > 0$, there exists some $T = T(\varepsilon, r) > 0$ so that $\varphi(s, t) < \varepsilon$ for all $s \leq r$ and $t \geq T$;
- for all $\varepsilon > 0$, there exists $\delta > 0$ so that $\varphi(s, t) \leq \varepsilon$ for all $s \leq \delta$ and all $t \geq 0$;

then there exists some $\beta \in \mathcal{KL}$ so that $\varphi(s, t) \leq \beta(s, t)$ for all $s \geq 0$, $t \geq 0$.

■

The following technical lemma will also be needed

**Lemma A.2** Suppose given a function $\delta \in \mathcal{K}_\infty$, a constant $C \geq 0$, and a function

$$T : (C, \infty) \times (0, \infty) \to \mathbb{R}_{\geq 0}$$

which satisfies the following two properties:

- $T$ is nondecreasing in its second argument, i.e., for all $\varepsilon > C$:

$$r_1 < r_2 \Rightarrow T(\varepsilon, r_1) \leq T(\varepsilon, r_2),$$

- for all $r > 0$, $\lim_{\varepsilon \searrow C} T(\varepsilon, r) = +\infty$.

Then, there is a function $\beta \in \mathcal{KL}$ with the following property: For each $r > 0$ and $t \geq 0$, there exists some

$$\varepsilon \in A_{r,t} := \{\varepsilon \in (C, \infty) \,|\, T(\varepsilon, r) \leq t\} \bigcup \{\infty\}, \tag{22}$$

such that

$$\min\{\varepsilon - 2C, \delta^{-1}(r)\} \leq \beta(r, t). \tag{23}$$

*Proof.* Introduce the following function:

$$\overline{\varphi}(r,t) := \inf A_{r,t}$$

defined for $r > 0$ and $t \geq 0$ (with $\overline{\varphi}(r,t) = \infty$ if $A_{r,t} = \{\infty\}$). The limit assumption on $T$ implies that there is, for each $r$ and $t$, some $\mu > 0$ such that $(C, C + \mu) \cap A_{r,t} = \emptyset$, from which it follows that $\overline{\varphi}(r,t) > C \geq 0$ for all $r, t$. Now let:

$$\varphi(r,t) := \min\{\overline{\varphi}(r,t) - C, \delta^{-1}(r)\}$$

for $r > 0$ and $t \geq 0$, and $\varphi(0, t) = 0$ for all $t \geq 0$, with the consistent convention that $\varphi(r,t) = \delta^{-1}(r)$ when $\overline{\varphi}(r,t) = \infty$. Note that, for each $r > 0$ and $\varepsilon > 0$, it holds that

$$[0 < s \leq r \quad \text{and} \quad t \geq T(\varepsilon + C, r)] \;\Rightarrow\; \varphi(s,t) \leq \varepsilon$$

(because $T(\varepsilon + C, s) \leq T(\varepsilon + C, r) \leq t \Rightarrow \varepsilon + C \in A_{s,t} \Rightarrow \overline{\varphi}(s,t) \leq \varepsilon + C$) and, for each $\varepsilon > 0$:

$$0 < s \leq \delta(\varepsilon) \;\Rightarrow\; \varphi(s,t) \leq \varepsilon$$

(because $\varphi(s,t) \leq \delta^{-1}(s) \leq \varepsilon$). Thus, by Proposition A.1, there is some $\overline{\beta} \in \mathcal{KL}$ so that $\varphi(r,t) \leq \overline{\beta}(r,t)$ for all $r \geq 0$ and all $t \geq 0$.

Define the new $\mathcal{KL}$ function $\beta := 2\overline{\beta}$. Pick any $r > 0$ and $t \geq 0$, and consider $\overline{\varphi}(r,t) \in (0, \infty]$. If $\overline{\varphi}(r,t) < \infty$, then there is (by definition of the function $\overline{\varphi}$) some $\varepsilon \in A_{r,t}$ such that $\varepsilon < 2\overline{\varphi}(r,t)$. From

$$\min\{\overline{\varphi}(r,t) - C, \delta^{-1}(r)\} \leq \overline{\beta}(r,t)$$

and $\varepsilon - 2C < 2\overline{\varphi}(r,t) - 2C$, one concludes (23). If, instead, $\overline{\varphi}(r,t) = \infty$, then $\delta^{-1}(r) \leq \overline{\beta}(r,t)$, and the conclusion holds with $\varepsilon = \infty$. ∎

*Proof.* (Proposition 2.6). The proof is a minor extension of that found in [19].

(i) : Suppose a set $S$ of trajectories satisfies the $(\mu, C)$-$\mathcal{KL}$-practical-IOS property with some $\mathcal{KL}$-function $\beta$. Define $\widehat{\beta}(s) := \beta(s, 0)$. Choose $\delta(\cdot) \in \mathcal{K}_\infty$ so that

$$\delta(\varepsilon) \leq \widehat{\beta}^{-1}(\varepsilon)$$

for $\varepsilon \in [0, \sup_{s \geq 0} \widehat{\beta}(s))$. Then for each $\varepsilon > 0$ and each trajectory $(\tau, u(\cdot), x(\cdot), y(\cdot)) \in S$, if $t_0 \in [0, \tau)$ is such that $x(t_0) \leq \delta(\varepsilon)$, then for all $t \in [t_0, \tau)$,

$$\begin{aligned} y(t) &\leq \max\{\beta(x(t_0), t - t_0), \|u\|_{[t_0, t]}, C\} \\ &\leq \max\{\beta(\delta(\varepsilon), 0), \|u\|_{[t_0, t]}, C\} \\ &\leq \max\{\varepsilon, \|u\|_{[t_0, t]}, C\}. \end{aligned}$$

Thus the set $S$ is $(\mu, C)$-uniformly practically stable. The $(\mu, C)$-uniform practical attractivity property follows immediately from $\beta(r, t) \to 0$ as $t \to \infty$.

(ii) : Suppose now the set $S$ satisfies the $(\mu, C)$-practical-IOS property. This means in particular that there exist a function $\delta \in \mathcal{K}_\infty$ and a function $T : (C, \infty) \times (0, \infty) \to \mathbb{R}_{>0}$ such that, for each trajectory in $S$, if we denote $r = x(t_0)$ and we pick any $t \in [t_0, \tau)$, the following two properties hold:

- $y(t) \leq \max\{\delta^{-1}(r), \|u\|_{[t_0, t]}, C\}$,

- $y(t) \leq \max\{\varepsilon, \|u\|_{[t_0, t]}\}$ for any $\varepsilon < \infty$ in the set $A_{r, t-t_0}$ in (22).

Without loss of generality, we may assume that $T$ is as in Lemma A.2. Indeed, we can always replace, if needed, $T$ by

$$\widehat{T}(\varepsilon, r) := \frac{r}{\varepsilon - C} + \inf\{T(\varepsilon', r') \mid r' \geq r, C < \varepsilon' \leq \varepsilon\}$$

which is even strictly increasing in $r$ and strictly decreasing in $\varepsilon$.

Pick $\beta$ as in Lemma A.2, any trajectory in $S$, and any $t \in [t_0, \tau)$. Denote $r = x(t_0)$. Pick $\varepsilon \in A_{r, t-t_0}$ as in the Lemma, for the given $r$ and $t$, so $\min\{\varepsilon - 2C, \delta^{-1}(r)\} \le \beta(r, t - t_0)$.

There are two cases to consider. If $\varepsilon - 2C < \delta^{-1}(r)$, necessarily $\varepsilon < \infty$ and so $y(t) \le \max\{\varepsilon, \|u\|_{[t_0, t]}\}$ gives
$$y(t) \le \max\{2C + \beta(r, t - t_0), \|u\|_{[t_0, t]}\}.$$
If, instead, $\varepsilon - 2C \ge \delta^{-1}(r)$, then $\delta^{-1}(r) \le \beta(r, t - t_0)$ gives that
$$y(t) \le \max\{\beta(r, t - t_0), \|u\|_{[t_0, t]}, C\}.$$
We conclude that
$$y(t) \le \max\{3\beta(r, t - t_0), \|u\|_{[t_0, t]}, 3C\}$$
for all $t \in [t_0, \tau)$. ∎